\documentclass[11pt]{article} 
\pdfoutput=1 
\usepackage{amssymb,amsmath,amsthm,latexsym,
ifthen, amsfonts,amsbsy,mathrsfs,bbm,mathtools}
\usepackage{graphicx}
\usepackage{color}
\usepackage{subfigure, latexsym,ifthen, amsfonts,amsbsy}
\usepackage{tikz}
\usetikzlibrary{arrows,shapes.arrows,calc,shapes.geometric,shapes.multipart,  
decorations.pathmorphing,positioning,swigs}

\RequirePackage[colorlinks]{hyperref}

 \usepackage{tablefootnote}
 \usepackage{float}
  \usepackage{url}

\usepackage[round,comma]{natbib}

\newtheorem{theorem}{Theorem}

\newtheorem{lemma}{Lemma}

\newcommand{\ind}{\mbox{$\perp \kern-5.5pt \perp$}}

\title{Assumptions and Bounds in the Instrumental Variable Model
}
\begin{document}

\author{Thomas S. Richardson\\
\em University of Washington
 \and James M. Robins\\
\em Harvard School of Public Health}
\maketitle

In this note we give proofs for results relating to  the  Instrumental Variable (IV) model with binary response $Y$ and binary treatment $X$,  but with an instrument $Z$ with $K$ states. These results were originally stated in {\it ACE Bounds; SEMS with Equilibrium Conditions} \citep{richardson:robins:imbens:2014}.

\section{Instrumental Variable Models}
We consider the Instrumental Variable (IV) model in which the treatment $X$ and response $Y$ are binary, taking values in $\{0,1\}$, while
 the instrument $Z$ takes values in $\{1,\ldots ,K\}$. We use the notation $X(z_k)$ to indicate
 $X(z\!=\!k)$, similarly $Y(x_i)$ for  $Y(x\!=\!i)$. 
 All of the models we consider will assume the individual level exclusion restriction
 $Y(x) = Y(x,z)$. We consider four different sets of independence assumptions:

\begin{itemize}\itemsep2pt
\item[(i)] $Z\; \ind\;  Y({x}_0),Y({x}_1),  X({z}_1),\ldots,X({z}_{K})$;
\item[(ii)] $Z\; \ind\;  Y({x}_0),Y({x}_1)$;
\item[(iii)] for $k \in \{1,\ldots ,K\}$, $i \in \{0,1\}$, $Z\; \ind\;  X(z_k),Y(x_i)$; 
\item[(iv)] there exists $U$ such that $U\; \ind\; Z$ and for $i \in \{0,1\}$, $Y(x_i)\; \ind\; X,Z \mid U$.
\end{itemize}

Condition (i) is joint independence of $Z$ and all the potential outcomes $Y(x_0)$, $Y(x_1)$, $X(z_1),\ldots ,X(z_K)$. 
This condition arises most naturally in contexts where the instrument is randomized via an external intervention.

Condition (ii) is implied by (i), but does not assume independence (or existence) of counterfactuals for $X$; see \citet{kitagawa:2021} for additional analysis of this model.

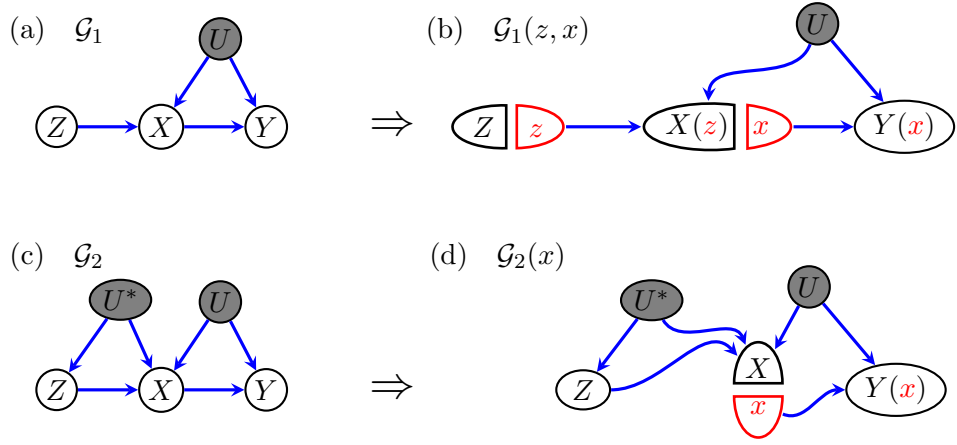
\begin{figure}[h]
\begin{tikzpicture}[>=stealth, ultra thick, node distance=2cm, inner sep=1.5,
    pre/.style={->,>=stealth,ultra thick,blue,line width = 1.2pt}]
\begin{scope}
\node[thick, name=z,shape=circle,style={draw}]
{$Z$
};
\node[above=1cm of z.center]{(a)\quad ${\cal G}_1$};
\node[thick,name=x,shape=circle,style={draw}, right=0.8cm of z] {
$X$
};
\node[thick, name=y,shape=circle,style={draw},right =0.8cm of x]
{$Y$
};
\draw[pre,->] (x) to (y);
\node[name=ca,right =0.4cm of x]{};
\node[thick,name=u,above=0.8cm of ca,shape=circle,style={draw},fill=gray] {
$U$
};
\draw[pre,->] (u) to (x);
\draw[pre,->] (u) to (y);
\draw[pre,->] (z) to (x);
\node[right =1cm of y,font=\LARGE]{$\Rightarrow$};
\end{scope}
\begin{scope}[xshift=6cm]
\node[line width = 1pt, inner sep=4pt, name=z, shape=swig vsplit, swig vsplit={line color right=red, gap=4pt}]{
         \nodepart{left}{$Z$}
         \nodepart{right}{\color{red}{$z$}}  
         };
\node[line width = 1pt,name=x, shape=swig vsplit,  right=1cm of z, swig vsplit={line color right=red, gap=5pt}]{
        \nodepart{left}{$X({\color{red}{z}})$}
        \nodepart{right}{$\color{red}{x}$}   };
\node[above=1cm of z.center]{(b)\quad ${\cal G}_1(z,x)$};
\node[thick, name=y,shape=ellipse,style={draw},right =0.8cm  of x]
{$Y({\color{red} {x}})$
};
\draw[pre,->] (z) to (x);
\node[name=cb,left =0.4cm of y]{};
\node[thick,name=u,above=1cm of cb,shape=circle,style={draw},fill=gray] {
$U$
};
\draw[pre,->] (u) to[out=250, in=90](x);
\draw[pre,->] (u) to (y);
\draw[pre,->] (x) to (y);
\end{scope}
\begin{scope}[yshift=-3.5cm]
\node[thick, name=z,shape=circle,style={draw}]
{$Z$
};
\node[above=1.5cm of z.center]{(c)\quad ${\cal G}_2$};
\node[thick,name=x,shape=circle,style={draw}, right=0.8cm of z] {
$X$
};
\node[thick, name=y,shape=circle,style={draw},right =0.8cm of x]
{$Y$
};
\draw[pre,->] (x) to (y);
\node[name=ca,right =0.4cm of x]{};
\node[thick,name=u,above=0.8cm of ca,shape=circle,style={draw},fill=gray] {
$U$
};
\node[thick,name=ustar,above right=0.8cm and 0.35cm of z,shape=ellipse,style={draw},fill=gray] {
$U^*$
};
\draw[pre,->] (ustar) to (x);
\draw[pre,->] (ustar) to (z);
\draw[pre,->] (u) to (x);
\draw[pre,->] (u) to (y);
\draw[pre,->] (z) to (x);
\node[right =1cm of y,font=\LARGE]{$\Rightarrow$};
\end{scope}
\begin{scope}[xshift=7cm,yshift=-3.5cm]
\node[thick,name=z,shape=ellipse, style={draw},text width=4mm,align=center] {
$Z$
};
\node[above left=1.5cm and 0.15cm of z.center]{(d)\quad ${\cal G}_2(x)$};
\node[line width = 1pt,name=x, shape=swig hsplit,  right=1.6cm of z, swig hsplit={line color lower=red, gap=5pt}]{
        \nodepart{upper}{$X$}
        \nodepart{lower}{$\color{red}{x}$}   };

\node[thick, name=y,shape=ellipse,style={draw},right =0.8cm  of x]
{$Y({\color{red} {x}})$
};
\draw[pre,->] (z) to[out=0,in=130] (x);
\node[name=cb,left =0.4cm of y]{};
\node[thick,name=u,above=1cm of cb,shape=circle,style={draw},fill=gray] {
$U$
};
\node[thick,name=ustar,above right=0.8cm and 0.35cm of z,shape=ellipse,style={draw},fill=gray] {
$U^*$
};
\draw[pre,->] (ustar) to[out=-60, in=110] (x);
\draw[pre,->] (ustar) to (z);
\draw[pre,->] (u) to (x);
\draw[pre,->] (u) to (y);
\draw[pre,->] (x) to[out=-30,in=180] (y);
\end{scope}
\end{tikzpicture}
\caption{{\rm (a)} IV model with no confounding between $Z$ and $X$; {\rm (b)} SWIG representing $P(Z, X(z),Y(x),U)$;
{\rm (c)} IV model with confounding between $Z$ and $X$; {\rm (d)} SWIG representing $P(Z, X,Y(x),U,U^*)$.\label{fig:swig} }
\end{figure}

Condition (iii) is also a subset of the independences in (i); note that none of the $2K$ independences in (iii) involve potential outcomes from different worlds.\footnote{In other words,
they do not involve both $Y(x_0)$ and $Y(x_1)$, nor $X(z_i)$ and $X(z_j$) for $i\neq j$.}
Note that the assumption (iii) may be read (via d-separation) from the Single-World Intervention Graph (SWIG)\footnote{See \citet{richardson:robins:2013} for details.}
${\cal G}_1(z,x)$, depicted in  Figure~\ref{fig:swig}(b), which represents the factorization of $P(Z,X(z),Y(x),U)$, implied by the 
Finest Fully Randomized Causally Interpreted Structured Tree Graph (FFRCISTG) associated with the DAG.

Lastly (iv) consists of only three independence statements, but does assume the existence of an unobserved variable $U$ that 
is sufficient to control for confounding between $X$ and $Y$. No assumption is made concerning the existence of counterfactuals
$X(z)$; confounding variables ($U^*$) between $Z$ and $X$ are permitted (so long as $U^* \ind U$). The DAG ${\cal G}_2$ and
corresponding SWIG ${\cal G}_2(x)$ are shown in Figure~\ref{fig:swig}(c),(d).

\section{Characterization of the IV Models}  
  
\begin{theorem}\label{thm:one} 
Under any of the assumptions {\rm (i)}, {\rm (ii)}, {\rm (iii)},  or {\rm (iv)} the joint distribution $P(Y(x_0), Y(x_1))$ is characterized by
the $4K$ pairs of inequalities:
\begin{align}
 P(Y(x_i)\!=\!y) &\leq  P(Y\!=\!y, X\!=\!i \,|\, Z\!=\! z) + P(X\!=\!1-i \,|\, Z\!=\!z),\label{eq:marg}\\
 P(Y(x_0)\!=\!y, Y(x_1)\!=\!\tilde{y}) &\leq P(Y\!=\!y, X\!=\!0\,|\, Z\!=\! z) +  P(Y\!=\!\tilde{y}, X\!=\!1\,|\, Z\!=\! z),\label{eq:joint}
\end{align}
with $y, \tilde{y}, i \in \{0,1\}$ and $z\in \{1,\ldots ,K\}$.
\end{theorem} 

Before proving this we first establish three useful Lemmas. 

For each of the models (i),\ldots,(iv), we consider the possible joint distributions 
\begin{equation}\label{eq:base-joint}
P(Z,X,Y(x_0),Y(x_1))
\end{equation}
implied by each model; we will use $M_1, \ldots, M_4$ to denote these. 

Note that although the models differ regarding the random variables that they posit, for example $M_2$ does not assume that $X(z)$ exists, they all contain the variables in  (\ref{eq:base-joint}); under $M_3$, $X=X(Z)$ by consistency. The observed distribution is then directly determined via consistency:
\[
P(Y\!=\!y, X\!=\!i \mid Z\!=\!z) = P(Y(x_i)\!=\!y, X\!=\!i \mid Z\!=\!z),
\]
where the RHS is a function of (\ref{eq:base-joint}).

\begin{lemma}\label{prop:contains}
The following containment relations hold:
\begin{align}
M_1 &\subseteq M_j; \quad j=2,3,4;\\
M_2 &\subseteq M_4.
\end{align}
\end{lemma}

\noindent{\it Proof:} $M_1 \subseteq M_2$ and $M_1 \subseteq M_3$ follow immediately by definition.
$M_2 \subseteq M_4$ follows by noting that under $M_2$ we may take $U= (Y(x_0),Y(x_1))$.
\hfill$\Box$

\bigskip

\begin{lemma} \label{prop:m4-z-ind}
$M_4$ implies $Z \ind Y(x_i)$, for $i \in \{0,1\}$.
\end{lemma}

\noindent{\it Proof:} $Y(x_i) \ind X,Z \mid U$ $\Rightarrow$ $Z \ind Y(x_i) \mid U$. Together with $Z \ind U$, this implies $Z\ind Y(x_i), U$, from which
the conclusion follows.\hfill$\Box$

\begin{lemma} \label{prop:marg-holds}
If $Z \ind Y(x_i)$ under $P$ then {\rm (\ref{eq:marg})} holds.
\end{lemma}

\noindent{\it Proof:} 
\begin{align*}
P(Y(x_i) = y) &= P(Y(x_i)\!=\!y \mid Z\!=\!z)\\
&= P(Y(x_i)\!=\!y, X=i \mid Z\!=\!z) + P(Y(x_i)\!=\!y, X=1-i \mid Z\!=\!z)\\
&= P(Y\!=\!y, X=i \mid Z\!=\!z) + P(Y(x_i)\!=\!y, X=1-i \mid Z\!=\!z)\\
&\leq  P(Y\!=\!y, X=i \mid Z\!=\!z) + P(X\!=\!1-i \mid Z\!=\!z).
\end{align*}
Here the first line follows from the independence, the second from the law of total probability, the third from consistency,
and the fourth via a probability inequality.
\hfill$\Box$

\bigskip

We define:
\begin{align}
\phi\;:\; P(Z,X,Y(x_0),Y(x_1)) \;\mapsto\; \left( P(Y(x_0),Y(x_1)), P(X,Y\mid Z)\right);
\end{align}
thus $\phi (P)$ gives the marginal distribution over the potential outcomes, and the observed distributions conditional on $Z$.
By extension $\phi(M_i)$ will denote the image, under $\phi$, of the set of distributions in $M_i$.
Let
${\cal T}$ denote the set of pairs 
 $ \left( P(Y(x_0),Y(x_1)), P(X,Y\mid Z)\right)$ that obey (\ref{eq:marg}) and (\ref{eq:joint}).

\bigskip

\noindent{\it Proof:} (Theorem \ref{thm:one})\par 
\noindent Theorem \ref{thm:one} is equivalent to the claim that $\phi(M_i)={\cal T}$, for $i=1,\ldots ,4$.

We prove this in three steps by establishing: (1) $\phi(M_3) \subseteq {\cal T}$; (2) $\phi(M_4) \subseteq {\cal T}$;  (3) ${\cal T} \subseteq \phi(M_1)$.
The conclusion then follows from Lemma \ref{prop:contains} since, by definition, $\phi(M_1) \subseteq \phi(M_j)$ for $j=2,3,4$.

\subsection*{(1) Proof that $\phi(M_3) \subseteq {\cal T}$}

Since $Z\ind Y(x_i)$ under $M_3$, (\ref{eq:marg}) holds by Lemma \ref{prop:marg-holds}. 
Thus it only remains to show (\ref{eq:joint}).

\begin{eqnarray*}
\lefteqn{P(X\!=\!0,Y\!=\!y\,|\, Z\!=\!z) +  P(X\!=\!1,Y\!=\!\tilde{y} \,|\, Z\!=\!z)}\\[6pt]
&=& P(X(z)\!=\!0,Y(x_0)\!=\!y \,|\, Z\!=\!z) \\[-2pt]
&&\quad\quad \quad\quad \quad\quad +  P(X(z)\!=\!1,Y(x_1)\!=\!\tilde{y} \, |\, Z\!=\!z)\\[6pt]
&=&   P(X(z)\!=\!0,Y(x_0)\!=\!y) +  P(X(z)\!=\!1,Y(x_1)\!=\!\tilde{y})\\[6pt]
&\geq & P(X(z)\!=\!0,Y(x_0)\!=\!y,Y(x_1)\!=\!\tilde{y}) +\\[-2pt]
&&\quad\quad\quad\quad\quad\quad\quad\quad  P(X(z)\!=\!1,Y(x_0)\!=\!y,Y(x_1)\!=\!\tilde{y})\\[6pt]
&= & P(Y(x_0)\!=\!y,Y(x_1)\!=\!\tilde{y}) 
\end{eqnarray*}
The first equation is by consistency; the second follows since  $M_3$ assumes $X(z),Y(x_i) \ind Z$; the third is a probability inequality.
\hfill$\Box$

\subsection*{(2) Proof that $\phi(M_4) \subseteq {\cal T}$}

By Lemma \ref{prop:m4-z-ind}, $Z\ind Y(x_i)$ holds under $M_4$, so (\ref{eq:marg}) holds by Lemma \ref{prop:marg-holds}. 
It remains to show (\ref{eq:joint}).

\begin{eqnarray*}
\lefteqn{ P(X\!=\!0,Y\!=\!y\,|\, Z\!=\!z) +  P(X\!=\!1,Y\!=\!\tilde{y} \,|\, Z\!=\!z)}\\[6pt]
&=& \sum_{u} P(Y\!=\!y, X\!=\!0, U\!=\!u\,|\, Z\!=\!z) + P(Y\!=\!\tilde{y}, X\!=\!1, U\!=\!u\,|\, Z\!=\!z)\\[6pt]
&=& \sum_{u}\left\{ P(Y\!=\!y\,|\, X\!=\!0, U\!=\!u, Z\!=\!z)P(X\!=\!0\,|\, U\!=\!u, Z\!=\!z) +\right.  \\[-8pt]
&&  \kern20pt \left. P(Y\!=\!\tilde{y}\,|\, X\!=\!1, U\!=\!u, Z\!=\!z)P(X\!=\!1\,|\, U\!=\!u, Z\!=\!z)\right\}P(U\!=\!u  \,|\, Z\!=\!z)\\[6pt]
&=& \sum_{u}\left\{ P(Y(x_0) \!=\!y\,|\, X\!=\!0, U\!=\!u, Z\!=\!z)P(X\!=\!0\,|\, U\!=\!u, Z\!=\!z) +\right.  \\[-8pt]
&&  \kern20pt \left. P(Y(x_1) \!=\!\tilde{y}\,|\, X\!=\!1, U\!=\!u, Z\!=\!z)P(X\!=\!1\,|\, U\!=\!u, Z\!=\!z)\right\}P(U\!=\!u  \,|\, Z\!=\!z)\\[6pt]
&=& \sum_{u}\left\{ P(Y(x_0) \!=\!y\,|\,  U\!=\!u)P(X\!=\!0\,|\, U\!=\!u, Z\!=\!z) +\right.  \\[-8pt]
&&  \kern20pt \left. P(Y(x_1) \!=\!\tilde{y}\,|\,  U\!=\!u)P(X\!=\!1\,|\, U\!=\!u, Z\!=\!z)\right\}P(U\!=\!u )\\[6pt]
&\geq &  \sum_{u}\left\{ P(Y(x_0)\!=\!y,Y(x_1)\!=\!\tilde{y}\,|\, U\!=\!u)P(X\!=\!0\,|\, U\!=\!u, Z\!=\!z) +\right.  \\[-8pt]
&& \kern20pt  \left. P(Y(x_0)\!=\!y,Y(x_1)\!=\!\tilde{y}\,|\, U\!=\!u)P(X\!=\!1\,|\, U\!=\!u, Z\!=\!z)\right\}P(U\!=\!u)\\[6pt]
&= &  \sum_{u} P(Y(x_0)\!=\!y,Y(x_1)\!=\!\tilde{y}\,|\, U\!=\!u)P(U\!=\!u)\\[6pt]
&= & P(Y(x_0)\!=\!y,Y(x_1)\!=\!\tilde{y})
\end{eqnarray*}

The first and second equations are algebra; the third follows by consistency,
the fourth since $U \ind Z$ and $Y(x) \ind X, Z \mid U$;
 the fifth is a probability inequality; the sixth is algebra.\hfill$\Box$

\subsection*{Proof that ${\cal T} \subseteq \phi(M_1)$}

We will show that given a pair of distributions 
\[
\widetilde{Q} \equiv \left( \widetilde{P}(Y(x_0),Y(x_1)), \widetilde{P}(X,Y\mid Z)\right) \in {\cal T}
\]
then there exists a joint distribution
$P(Z,X(z_1),\ldots ,X(z_k),Y(x_0),Y(x_1))$ in model $M_1$, and such that 
$\phi(P(Z,X,Y(x_0),Y(x_1)))  = \widetilde{Q}$.

\medskip

We first note the following construction:

\begin{lemma}\label{lem:combine-k}
Given a set of $K$ distributions $P_k(X(z_k),Y(x_0),Y(x_1))$, for $k \in \{1,\ldots ,K\}$ that agree on the common marginals
so that $P_k(Y(x_0),Y(x_1))=P_{k^*}(Y(x_0),Y(x_1))$ for all $k,k^*$, then there exists a single joint distribution:
\[
P(X(z_1),\ldots ,X(z_k),Y(x_0),Y(x_1))
\]
that agrees with each of these $K$ marginals so that for all $k$,
\[
P_k(X(z_k),Y(x_0),Y(x_1)) =  P(X(z_k),Y(x_0),Y(x_1))).
\]
\end{lemma}

\noindent{\it Proof}: 

We may form a joint distribution
\begin{align*}
P(X(z_1),\ldots, X(z_k), Y(x_0),Y(x_1)) &= \frac{\prod_{k=1}^K P_{k}(X(z_k),Y(x_0), Y(x_1))}
{P_{1}(Y(x_0),Y(x_1))^{K-1}}.
\end{align*}
The resulting distribution $P^*$ agrees with each $P_k$ on the $(X(z_k),Y(x_0),Y(x_1))$
margin.\footnote{Though not important for our argument we note that $P^*$ enforces the joint conditional independence
of the $X(z)$ counterfactuals given $Y(x_0),Y(x_1)$:
\[
X(z_1)\;\ind\; X(z_2)\;\ind\; \cdots \;\ind\; X(z_k) \;\mid\; \{Y(x_0),Y(x_1)\}.
\]}
\hfill$\Box$
\bigskip

\noindent Under the IV model (i), we have the following equality:
\begin{align}
P(Y(x_j)=i, X(z_k)=j) &= P(Y(x_j)=i, X(z_k)=j \mid Z=k)\nonumber\\
&= P(Y=i, X=j \mid Z=k).\label{eq:z-obs-constraint}
\end{align}
Lemma \ref{lem:combine-k} implies that it is sufficient to consider separately each level of $Z=k \in \{1,\ldots, K\}$:
If we can construct $K$ marginal distributions over $(X(z_k),Y(x_0),Y(x_1))$ each of which obeys the corresponding constraints (\ref{eq:z-obs-constraint})
and that agree on the $(Y(x_0),Y(x_1))$ margin, 
then these may be combined to form a single joint distribution.

It thus remains to show that, given a pair $\widetilde{P}(Y(x_0),Y(x_1))$, $\widetilde{P}(X,Y\mid Z\!=\!k)$, obeying 
 (\ref{eq:marg}) and (\ref{eq:joint}), these inequalities are sufficient to ensure that there exists a joint distribution 
 $P(X(z_k),Y(x_0),Y(x_1)))$, such that:

\begin{align}
P(Y(x_0)\!=\!y, Y(x_1)\!=\!y^*) &= \widetilde{P}(Y(x_0)\!=\!y, Y(x_1)\!=\!y^*),\label{eq:linkp-to-ptilde}\\
P(X(z_k)\!=\!j, Y(x_j)\!=\!y) &= \widetilde{P}(Y\!=\!y, X\!=\!j \mid Z\!=\!k).\label{eq:linkz-to-ptilde}
\end{align}

The existence of such a distribution could be established by appeal to Strassen's Theorem.
However, here we opt to do this directly by elementary arguments.

We first claim that  the following probabilities are sufficient to parametrize $P(X(z_k),Y(x_0),Y(x_1)))$:
\begin{align}
P(Y(x_0)\!=\!y, Y(x_1)\!=\!y^*),&\quad y,y^* \in \{0,1\},\label{eq:givens1}\\ 
P(X(z_k)=j,Y(x_j)\!=\!y), &\quad j,y \in \{0,1\},\label{eq:givens2}\\
 \color{blue}{P(X(z_k)=0, Y(x_0)=0, Y(x_1)=0)}.&\label{eq:givens3}
\end{align}
By hypothesis, the quantities in (\ref{eq:givens1}) and (\ref{eq:givens2}), are specified by (\ref{eq:linkp-to-ptilde}) and (\ref{eq:linkz-to-ptilde}) respectively.
To establish the claim we express each of the $8$ joint probabilities in terms of the quantities (\ref{eq:givens1}), (\ref{eq:givens2}) and (\ref{eq:givens3}):
\begin{align}
\MoveEqLeft
{P(X(z_k)=0,Y(x_0)=0, Y(x_1)=0) =}\label{eq:001}\\
&= \color{blue}{P(X(z_k)=0,Y(x_0)=0, Y(x_1)=0)}
\nonumber\\[4pt]
\MoveEqLeft
{P(X(z_k)=0,Y(x_0)=0, Y(x_1)=1) =}\label{eq:011}\\
&= P(X(z_k)=0,Y(x_0)=0) - \color{blue}{P(X(z_k)=0,Y(x_0)=0, Y(x_1)=0)}
\nonumber\\[4pt]
\MoveEqLeft
{P(X(z_k)=0,Y(x_0)=1, Y(x_1)=0) =}\label{eq:101}\\
&= P(Y(x_1)=0) - P(X(z_k)=1,Y(x_1)=0) - \color{blue}{P(X(z_k)=0,Y(x_0)=0, Y(x_1)=0)}\nonumber\
\\[4pt]
\MoveEqLeft
{P(X(z_k)=0,Y(x_0)=1, Y(x_1)=1) =}\label{eq:111}\\
&= P(X(z_k)=0,Y(x_0)=1) + P(X(z_k)=1,Y(x_1)=0)\nonumber\\
& \kern120pt - P(Y(x_1)=0) + \color{blue}{P(X(z_k)=0,Y(x_0)=0, Y(x_1)=0)}\nonumber\\[4pt]
\MoveEqLeft
{P(X(z_k)=1,Y(x_0)=0, Y(x_1)=0) =}\label{eq:002}\\
&= P(Y(x_0)=0, Y(x_1)=0) - \color{blue}{P(X(z_k)=0,Y(x_0)=0, Y(x_1)=0)}
\nonumber\\[4pt]
\MoveEqLeft
{P(X(z_k)=1,Y(x_0)=0, Y(x_1)=1) =\label{eq:012}}\\
&= P(Y(x_0)=0, Y(x_1)=1) - P(X(z_k)=0,Y(x_0)=0)\nonumber\\ 
&\kern120pt+ \color{blue}{P(X(z_k)=0,Y(x_0)=0, Y(x_1)=0)}\nonumber\
\\[4pt]
\MoveEqLeft
{P(X(z_k)=1,Y(x_0)=1, Y(x_1)=0) =}\label{eq:102}\\
&=  P(X(z_k)=1,Y(x_1)=0) - P(Y(x_0)=0, Y(x_1)=0)\nonumber\\
& \kern120pt+ \color{blue}{P(X(z_k)=0,Y(x_0)=0, Y(x_1)=0)}\nonumber\
\\[4pt]
\MoveEqLeft
{P(X(z_k)=1,Y(x_0)=1, Y(x_1)=1) =}\label{eq:112}\\
&= P(X(z_k)=0,Y(x_0)=0)+ P(X(z_k)=1,Y(x_1)=1)\nonumber \\
& \kern60pt - P(Y(x_0)=0, Y(x_1)=1) - \color{blue}{P(X(z_k)=0,Y(x_0)=0, Y(x_1)=0)}\nonumber
\end{align}

Note that $P(X(z_k)=0,Y(x_0)=0, Y(x_1)=0)$ is present in every one of (\ref{eq:001})--(\ref{eq:112}).

Further note that the set of valid distributions $P(X(z_k),Y(x_0), Y(x_1))$ is determined by the constraints:
\begin{align}
P(X(z_k),Y(x_0), Y(x_1)) &\geq 0\label{eq:positivity} \\
\sum_{j,y,\tilde{y}} P(X(z_k)=j,Y(x_0)=y, Y(x_1)=\tilde{y}) &= 1. \label{eq:sum-to-one}
\end{align}

Given any values for the quantities in  (\ref{eq:givens1}), (\ref{eq:givens2}), (\ref{eq:givens3}), if we compute $ P(X(z_k)=i,Y(x_0)=y, Y(x_1)=\tilde{y})$ using
equations (\ref{eq:001})--(\ref{eq:112}), it is simple to see that these will obey (\ref{eq:sum-to-one}), since the sum of these terms is equal to the sum of the four terms
in (\ref{eq:givens2}).

Thus it merely remains to show that there exists a possible value for $P(X(z_k)=0,Y(x_0)=0, Y(x_1)=0)$  compatible with the eight positivity constraints (\ref{eq:positivity}).
Re-arranging each of these inequalities reveals that we have four lower bounds, arising from (\ref{eq:001}),  (\ref{eq:111}), (\ref{eq:012}) and (\ref{eq:102}):
\begin{align}
\MoveEqLeft{\color{blue}{P(X(z_k)=0,Y(x_0)=0, Y(x_1)=0)}}\label{eq:lower-on-p000}\\
& \geq 0\nonumber \\
&\geq    P(Y(x_1)=0) - P(X(z_k)=0,Y(x_0)=1) - P(X(z_k)=1,Y(x_1)=0)\nonumber\\
&\geq  P(X(z_k)=0,Y(x_0)=0) - P(Y(x_0)=0, Y(x_1)=1)\nonumber \\
&\geq   P(Y(x_0)=0, Y(x_1)=0) - P(X(z_k)=1,Y(x_1)=0)\nonumber
\end{align}
and also four upper bounds from (\ref{eq:011}), (\ref{eq:101}), (\ref{eq:002}) 
and (\ref{eq:112}):
\begin{align}
\MoveEqLeft{\color{blue}{P(X(z_k)=0,Y(x_0)=0, Y(x_1)=0)}}\label{eq:upper-on-p000}\\
&\leq  P(X(z_k)=0,Y(x_0)=0) \nonumber\\
&\leq   P(Y(x_1)=0) - P(X(z_k)=1,Y(x_1)=0) \nonumber\\
& \leq P(Y(x_0)=0, Y(x_1)=0)\nonumber\\
&\leq    P(X(z_k)=0,Y(x_0)=0)+ P(X(z_k)=1,Y(x_1)=1) - P(Y(x_0)=0, Y(x_1)=1).\nonumber
\end{align}

Combining the upper and lower bounds now results in sixteen inequalities.
There are four each of the form $P(Y(x_j)=y, X(z_k)=j)\geq 0$ and $P(Y(x_0)=y, Y(x_0)=\tilde{y})\geq 0$,
which are automatically satisfied by the respective distributions; see Appendix \ref{app:one} for the detailed calculation. The remaining 8 inequalities are equivalent to 
 (\ref{eq:marg}) and (\ref{eq:joint}) for fixed $Z=k$, using (\ref{eq:linkz-to-ptilde}).
 Thus  it follows that if 
 \[
 \left(\tilde{P}(Y(x_0),Y(x_1)), \tilde{P}(X,Y\mid Z)\right) \in {\cal T},
 \]
then  there is always a non-empty set of possible values for $P(X(z_k)=0,Y(x_0)=0, Y(x_1)=0)$ compatible with 
\[
\left(\tilde{P}(Y(x_0),Y(x_1)), \tilde{P}(X,Y\mid Z=z)\right).
\]
 Consequently, we can always construct a trivariate joint $P(X(z_k),Y(x_0),Y(x_1)))$ compatible with $(\tilde{P}(Y(x_0),Y(x_1)), \tilde{P}(X,Y\mid Z=z) $ via  (\ref{eq:linkp-to-ptilde}) and (\ref{eq:linkz-to-ptilde}). 

 Since, by construction, these marginal distributions agree on $\widetilde{P}(Y(x_0),Y(x_1))$, by Lemma \ref{lem:combine-k}, these trivariate joints may be combined into a single joint distribution
\[
P(X(z_1),\ldots, X(z_K),Y(x_0),Y(x_1))).
\]
Multiplying this by an arbitrary independent marginal distribution for $P(Z)$ then creates a distribution obeying (i)
that satisfies (\ref{eq:linkp-to-ptilde}) and (\ref{eq:linkz-to-ptilde}) with respect to pair $\left( \widetilde{P}(Y(x_0),Y(x_1)), \widetilde{P}(X,Y\mid Z\!=\!k)\right)$.
This establishes that ${\cal T} \subseteq \phi(M_1)$ as required. \hfill$\Box$

\eject

\section{Bounds on  $P(Y(x_0))$, $P(Y(x_1))$ and $ACE(X\!\rightarrow\! Y)$}

We next provide bounds on  $P(Y(x_0))$ and $P(Y(x_1))$ and $ACE(X\!\rightarrow\! Y)\equiv P(Y(x_1)=1) - P(Y(x_0)=1)$.

\begin{theorem}\label{thm:two}
 Under any of the assumptions {\rm (i)}, {\rm (ii)}, {\rm (iii)},  {\rm (iv)} for all $i,j \in \{0,1\}$, $P(Y(x_i) = j) \leq g(i,j)$, where:
\begin{align}
g(i,j) &\equiv \min\left\{ \min_{z} \left[\vphantom{\hat{P}} P(X\!=\!i, Y\!=\!j \,|\, Z\!=\!z) + P(X\!=\!1-i \,|\, Z\!=\!z)\right], \right.\label{eq:gijbound}\\[3pt]
&\kern15pt\quad \min_{z, \tilde{z}:\, z \neq \tilde{z}}  \left[\vphantom{\hat{P}} P(X\!=\!i, Y\!=\!j \,|\, Z\!=\!z) +  P(X\!=\!1-i, Y\!=\!0 \,|\, Z\!=\!z)\right. \nonumber \\[-3pt]
 &\quad\quad\kern60pt  +\;\left.\vphantom{ \min_{z} }\left.\vphantom{\hat{P}}  P(X\!=\!i, Y\!=\!j \,|\, Z\!=\!\tilde{z}) +  P(X\!=\!1-i, Y\!=\!1 \,|\, Z\!=\!\tilde{z})\right] \right\}. \nonumber
\end{align}
Furthermore, $P(Y(x_0))$ and $P(Y(x_1))$ are variation independent. Consequently:
\begin{equation}
1-g(1,0)-g(0,1) \;\leq\; ACE(X\!\rightarrow\! Y) \;\leq\; g(0,0)+g(1,1)-1.\label{eq:ace-bnds}
\end{equation}
These bounds are sharp.
\end{theorem} 

\begin{table}
 \begin{center}
 \begin{tabular}{cc|c}
 \multicolumn{2}{c|}{Potential Outcomes} & Type \\
   $Y(x_0)$ & $Y(x _1)$ & \\
 \hline
  $0$ & $0$ & Never Recover (NR)\\
    $0$ & $1$ & Helped (HE)\\
  $1$ & $0$ & Hurt (HU)\\
  $1$ & $1$ & Always Recover (AR)
  \end{tabular}
  \end{center}
\caption{Table of potential response types. \label{tab:types}}
\end{table}

%

\noindent{\it Proof:} (Theorem \ref{thm:two})\par
The set of joint distributions $P(Y(x_0),Y(x_1))$ is characterized by the three quantities:
\[
\pi_0\equiv P(Y(x_0)=1),\quad \pi_1\equiv P(Y(x_1)=1),\quad  \pi_{AR}\equiv P(Y(x_0)=1, Y(x_1)=1).
\]
Here AR denotes {\it `Always Recover'}, see Table \ref{tab:types}. 
To find the set of possible marginal probabilities $(\pi_0,\pi_1)$ we simply re-write the inequalities
(\ref{eq:marg}) and (\ref{eq:joint}) in terms of $\pi_0$, $\pi_1$ and $\pi_{AR}$ and then eliminate
$\pi_{AR}$. Substituting for the joint probabilities $P(Y(x_0),Y(x_1))$ in (\ref{eq:joint}), together
with the constraint that probabilities are non-negative leads to: 
\begin{eqnarray*}
0\;\;\leq & 1-\pi_0-\pi_1+\pi_{AR} &\leq\;\; P(y_0\mid z_i);\\
0\;\;\leq & \pi_0 - \pi_{AR} &\leq\;\; P(x_0,y_1\mid z_i) + P(x_1,y_0 \mid z_i);\\
0\;\; \leq& \pi_1 - \pi_{AR} & \leq\;\;  P(x_1,y_1 \mid z_i) + P(x_0,y_0\mid z_i);\\
0\;\;\leq& \pi_{AR} &\leq\;\; P(y_1 \mid z_i).
\end{eqnarray*}
In addition, we have the inequalities from (\ref{eq:marg}) which may be re-written as:
\begin{eqnarray}
P(x_0,y_1 \mid z_i) \;\;\leq &\pi_0& \leq\;\; 1- P(x_0,y_0\mid z_i);\label{eq:natbounds1}\\
P(x_1,y_1\mid z_i) \;\;\leq &\pi_1& \leq\;\; 1-P(x_1,y_0\mid z_i).\label{eq:natbounds2}
\end{eqnarray}
Re-arranging the former set leads to the following upper and lower bounds on $\pi_{AR}$:
\begin{equation*}
\begin{array}{crclc}
(\hbox{nr}) & \pi_0+\pi_1 - 1&\leq  \pi_{AR} \leq &  \pi_0+\pi_1 + P(y_0\mid z_i) - 1; & (\hbox{NR}_i)\\
(\hbox{hu}_i) &\pi_0- P(x_0,y_1\mid z_i) - P(x_1,y_0 \mid z_i) &\leq   \pi_{AR} \leq & \pi_0; & (\hbox{HU})\\
(\hbox{he}_i)&\pi_1  - P(x_1,y_1 \mid z_i) - P(x_0,y_0\mid z_i)  & \leq  \pi_{AR}  \leq & \pi_1; & (\hbox{HE})\\
(\hbox{ar}) & 0&\leq \pi_{AR} \leq& P(y_1 \mid z_i).& (\hbox{AR}_i)
\end{array}
\end{equation*}

To make the argument simpler to follow we have labelled the upper and lower bounds, first to indicate the values of $(y_0,y_1)$ that are used with 
(\ref{eq:joint}) using the response `types' in Table \ref{tab:types}, and second whether we obtain different bounds for different levels of $Z$.

We may now eliminate $\pi_{AR}$ by requiring that in the last display each lower bound be less than or equal to each upper bound; a procedure that is known as Fourier-Motzkin elimination.
However,  in fact, we are only interested in those inequalities that are not already implied by (\ref{eq:natbounds1})
and  (\ref{eq:natbounds2}). 
The following four inequalities:
\[
(\hbox{nr}) \leq (\hbox{NR}_i),\quad (\hbox{hu}_i) \leq  (\hbox{HU}), \quad (\hbox{he}_i) \leq (\hbox{HE}),
\quad (\hbox{ar})  \leq (\hbox{AR}_i),
\]
trivially hold since $P(x,y\mid z) \geq 0$. Further, since $0 \leq \pi_j \leq 1$ for $j\in \{0,1\}$, the following four also hold trivially:
\[
(\hbox{nr}) \leq (\hbox{HU}),\quad (\hbox{nr}) \leq  (\hbox{HE}),\quad\quad
(\hbox{ar}) \leq (\hbox{HU}),\quad (\hbox{ar}) \leq  (\hbox{HE}).
\] 
This leaves eight comparisons, which we group according to the left-hand side:

\noindent One with the $(\hbox{nr})$ lower bound on $\pi_{AR}$:\par
\noindent $(\hbox{nr}) < (\hbox{AR}_i)$:
\begin{equation}
 \pi_0+\pi_1 - 1 \leq P(y_1 \mid z_i).\label{eq:nrAR}
 \end{equation}
 Three with the $(\hbox{hu}_i)$  lower bound on $\pi_{AR}$:\par
 \noindent $(\hbox{hu}_i) < (\hbox{NR}_{i^*})$:
\begin{equation}
 \pi_0- P(x_0,y_1\mid z_i) - P(x_1,y_0 \mid z_i) \leq \pi_0+\pi_1 + P(y_0\mid z_{i^*}) - 1;\label{eq:huiNRistarfirst}
\end{equation}
$(\hbox{hu}_i) < (\hbox{HE})$:
\begin{equation}
 \pi_0- P(x_0,y_1\mid z_i) - P(x_1,y_0 \mid z_i) \leq \pi_1;\label{eq:huiHE}
\end{equation}
$(\hbox{hu}_i) < (\hbox{AR}_{i^*})$:
\begin{equation}
 \pi_0- P(x_0,y_1\mid z_i) - P(x_1,y_0 \mid z_i) \leq P(y_1 \mid z_{i^*}).\label{eq:huiARistarfirst}
\end{equation}
 Three with the $(\hbox{he}_i)$  lower bound on $\pi_{AR}$:\par
 \noindent 
$(\hbox{he}_i) < (\hbox{NR}_{i^*})$:
\begin{equation}
\pi_1  - P(x_1,y_1 \mid z_i) - P(x_0,y_0\mid z_i) \leq \pi_0+\pi_1 + P(y_0\mid z_{i^*}) - 1; \label{eq:heiNRistarfirst}
\end{equation}
$(\hbox{he}_i) < (\hbox{HU})$:
\begin{equation}
\pi_1  - P(x_1,y_1 \mid z_i) - P(x_0,y_0\mid z_i) \leq \pi_0; \label{eq:heiHU}
\end{equation}
$(\hbox{he}_i) < (\hbox{AR}_{i^*})$:
\begin{equation}
\pi_1  - P(x_1,y_1 \mid z_i) - P(x_0,y_0\mid z_i) \leq P(y_1\mid z_{i^*}); \label{eq:heiARistarfirst}
\end{equation}
Finally, one with the $(\hbox{ar})$ lower bound ($=0$):\par
$(\hbox{ar}) < (\hbox{NR}_i)$:
\begin{equation}
0 \leq  \pi_0+\pi_1 + P(y_0\mid z_i) - 1.\label{eq:arNR}
\end{equation}
We will show that {\em all} of the constraints here, {\em other} than those involving two distinct levels of $Z$ --- namely, 
(\ref{eq:huiNRistarfirst}),  (\ref{eq:huiARistarfirst}), (\ref{eq:heiNRistarfirst}), (\ref{eq:heiARistarfirst}), with $i\neq i^\star$ ---
are implied by (\ref{eq:natbounds1}) and (\ref{eq:natbounds2}).

We first demonstrate that the four inequalities involving $i$ and $i^*$ are already
implied by (\ref{eq:natbounds1}) and (\ref{eq:natbounds2}) when $i=i^*$.
We may rewrite (\ref{eq:huiNRistarfirst}) and (\ref{eq:heiNRistarfirst}), respectively, as:
\begin{eqnarray}
P(y_1\mid z_{i^*}) - P(x_0,y_1\mid z_i) - P(x_1,y_0 \mid z_i) &\leq& \pi_1;\label{eq:huiNRistar}\\
P(y_1\mid z_{i^*}) - P(x_1,y_1 \mid z_i) - P(x_0,y_0\mid z_i) &\leq& \pi_0.\label{eq:heiNRistar}
\end{eqnarray}
In the case where $i=i^*$, these lower bounds reduce to:
\begin{eqnarray*}
P(x_1,y_1\mid z_i) - P(x_1,y_0 \mid z_i) &\leq& \pi_1,\\
P(x_0,y_1 \mid z_i) - P(x_0,y_0\mid z_i) &\leq& \pi_0,
\end{eqnarray*}
which are clearly implied by (\ref{eq:natbounds1}) and (\ref{eq:natbounds2}). Similarly, we may rewrite 
(\ref{eq:huiARistarfirst}) and (\ref{eq:heiARistarfirst}) as:
\begin{eqnarray}
 \pi_0 &\leq& P(y_1 \mid z_{i^*}) + P(x_0,y_1\mid z_i)+ P(x_1,y_0 \mid z_i),\label{eq:huiARistar}\\
 \pi_1   &\leq& P(y_1\mid z_{i^*}) + P(x_1,y_1 \mid z_i) + P(x_0,y_0\mid z_i).\label{eq:heiARistar}
\end{eqnarray}
In the case where $i=i^*$, it is easy to see that these upper bounds are implied by  (\ref{eq:natbounds1}) and (\ref{eq:natbounds2}).
Thus, if $Z$ has $K$ values then we obtain $K(K-1)$ upper and lower bounds on $\pi_0$ and $\pi_1$.
\bigskip

\noindent We now show that the four remaining inequalities --- namely, 
(\ref{eq:nrAR}), (\ref{eq:huiHE}), (\ref{eq:heiHU}) and (\ref{eq:arNR}) ---
 are implied by (\ref{eq:natbounds1}) and (\ref{eq:natbounds2}).

First, consider (\ref{eq:nrAR}). By the upper bounds in (\ref{eq:natbounds1}) and (\ref{eq:natbounds2}),
we have that:
\[
\pi_0 + \pi_1 -1 \leq 1- P(x_0,y_0\mid z_i) + 1- P(x_1,y_0\mid z_i) -1  = 1 - P(y_0 \mid z_i).
\]
Similarly, (\ref{eq:arNR}) follows directly from the lower bounds in (\ref{eq:natbounds1}) and (\ref{eq:natbounds2}).
Next consider (\ref{eq:huiHE}), which may be re-written as:
\[
 \pi_0- \pi_1  \leq P(x_0,y_1\mid z_i) + P(x_1,y_0 \mid z_i).
\]
Using the upper bound on $\pi_0$ from (\ref{eq:natbounds1}) and the lower bound on $\pi_1$
given by (\ref{eq:natbounds2}) we obtain:
\[
\pi_0 - \pi_1 \leq 1- P(x_0,y_0\mid z_{i}) - P(x_1,y_1\mid z_i) = P(x_0,y_1\mid z_i) + P(x_1,y_0 \mid z_i).
\]
Finally, we consider (\ref{eq:heiHU}), which can be written as:
\[
 \pi_1- \pi_0  \leq P(x_0,y_0 \mid z_i) + P(x_1,y_1 \mid z_i).
\]
This is implied by the upper bound on $\pi_1$ from (\ref{eq:natbounds2}) and the lower bound on $\pi_0$
from (\ref{eq:natbounds1}).

Summarizing, we have obtained the following bounds:
\begin{eqnarray*}
\pi_0 &\geq& \max_{i,i^*: i\neq i^*} \left\{ 
 P(x_0,y_1\mid z_i),\; 
P(y_1\mid z_{i^*}) - P(x_1,y_1 \mid z_i) - P(x_0,y_0\mid z_i)\right\},\\
\pi_0 &\leq& \min_{i,i^*: i\neq i^*} \left\{ 
 1- P(x_0,y_0\mid z_i),\; P(y_1 \mid z_{i^*}) + P(x_0,y_1\mid z_i)+ P(x_1,y_0 \mid z_i) \right\},
\end{eqnarray*}
\begin{eqnarray*}
\pi_1 &\geq& \max_{i,i^*: i\neq i^*} \left\{ 
P(x_1,y_1\mid z_i),\; 
P(y_1\mid z_{i^*}) - P(x_0,y_1\mid z_i) - P(x_1,y_0 \mid z_i)\right\},\\
\pi_1 &\leq& \min_{i,i^*: i\neq i^*} \left\{ 
1-P(x_1,y_0\mid z_i),\; 
P(y_1\mid z_{i^*}) + P(x_1,y_1 \mid z_i) + P(x_0,y_0\mid z_i) \right\},
\end{eqnarray*}
which can be written more compactly as (\ref{eq:gijbound}).

It is  an important and somewhat surprising  aspect of this problem that, although elimination of $\pi_{AT}$ produced inequalities involving constraints on $\pi_0 + \pi_1$ and $\pi_0 - \pi_1$, all four of these inequalities
were already implied by the original inequalities (\ref{eq:natbounds1}) and (\ref{eq:natbounds2}).
 Consequently, $\pi_0$ and $\pi_1$ are variation independent.\footnote{This appears to have been first noted by \citet{dawid:hsss} in the context of the model in which the instrument $Z$ is also binary.} This shows that sharp upper (lower) bounds on the ACE$(X\rightarrow Y) = \pi_1 - \pi_0$, may be obtained by considering the difference between the upper (lower) bound on $\pi_1$ and the lower (upper) bound on $\pi_0$, resulting in (\ref{eq:ace-bnds}). \hfill$\Box$
 \bigskip


In the special case where $Z$ is binary, these are the same bounds derived by \citet{balke1997bounds} under the stronger assumption that
$X$-counterfactuals existed and were independent of $Z$.
These bounds were also derived by Dawid, under a model that allowed for stochastic counterfactuals,
for both $X$ and $Y$, again assumed independent of $Z$ (though still restricted to binary $Z$).

Note that Pearl and Dawid express the bounds on $\pi_0$ and $\pi_1$ (for binary $Z$) in terms of four quantities; it is not immediately obvious that 
they may be re-expressed in the form given here, but this is the case.
Finally note that the restriction in the max to $i\neq i^*$ is purely for computational efficiency since from the previous analysis we know that 
the expression involving $i$ and $i^*$ can never be the upper or lower bound when  $i=i^*$.

\subsection{Bounds on the ACE$(X\rightarrow Y)$, when $X$ is binary}

It follows from the variation independence of $\pi_0$ and $\pi_1$ that bounds on the ACE $\pi_1 - \pi_0$ can easily be computed from the bounds on $\pi_0$ and $\pi_1$:
\[
l(\pi_1) - u(\pi_0)\leq  \pi_1 - \pi_0 \leq u(\pi_1)  - l(\pi_0)
\]
where $l(\cdot)$ and $u(\cdot)$ represent lower and upper bounds on the respective quantities under the observed distribution.

When $Z$ takes $K$ states we see little point in deriving more explicit expressions for the upper and lower bounds on the ACE\footnote{That is to say, expressions that are more explicit than (\ref{eq:gijbound}) and (\ref{eq:ace-bnds}).} since it is not clear that much is to be gained either by way of understanding or computational efficiency. However, we address one puzzle. In the case where all variables are binary, the upper and lower bounds on $\pi_0$ and $\pi_1$ each contain $4$ terms. It would thus be natural to expect that the upper and lower bounds on the ACE would each contain $4\times 4 = 16$ terms. However,
the bounds on the ACE reported by Pearl, and later re-derived by Dawid, contain only $8$ terms.  In Appendix \ref{app:two} we confirm by direct calculation that the additional $8$ inequalities are redundant (for both the upper and lower bounds). In addition, we relate the bounds on ACE$(X\rightarrow Y)$ given here by (\ref{eq:ace-bnds}) to those reported previously by \citet{pearl:2000, pearl:2009} and \citet{swanson:2018}.


\eject

\subsubsection*{Funding Information} 
The authors were supported by US Office of Naval Research grant N000141912446 and US National Institutes of Health Grant R01 AI032475; 
Richardson was also supported by the US National Science Foundation Grant CNS-0855230.


\bibliographystyle{chicago}
\bibliography{causal}

\appendix

\section{Detailed derivation of inequalities relating $P(Y(x_0),Y(x_1))$ and $P(X,Y\mid Z=k)$ \label{app:one}}

We here show explicitly how, for fixed $Z=k$, the four lower bounds given in (\ref{eq:lower-on-p000}),
\begin{align*}
\MoveEqLeft{\color{blue}{P(X(z_k)=0,Y(x_0)=0, Y(x_1)=0)}}\nonumber\\
\text{(a)} \quad & \geq 0 \\
\text{(b)} \quad &\geq    P(Y(x_1)=0) - P(X(z_k)=0,Y(x_0)=1) - P(X(z_k)=1,Y(x_1)=0)\\
\text{(c)} \quad &\geq  P(X(z_k)=0,Y(x_0)=0) - P(Y(x_0)=0, Y(x_1)=1)\\
\text{(d)} \quad &\geq   P(Y(x_0)=0, Y(x_1)=0) - P(X(z_k)=1,Y(x_1)=0),
\end{align*}
and  four upper bounds in (\ref{eq:upper-on-p000})
and (\ref{eq:112}),
\begin{align*}
\MoveEqLeft{\color{blue}{P(X(z_k)=0,Y(x_0)=0, Y(x_1)=0)}}\\
\text{(A)}\quad  &\leq  P(X(z_k)=0,Y(x_0)=0) \\
\text{(B)}\quad  &\leq   P(Y(x_1)=0) - P(X(z_k)=1,Y(x_1)=0) \\
\text{(C)}\quad  & \leq P(Y(x_0)=0, Y(x_1)=0)\\
\text{(D)}\quad  &\leq    P(X(z_k)=0,Y(x_0)=0)+ P(X(z_k)=1,Y(x_1)=1) - P(Y(x_0)=0, Y(x_1)=1),
\end{align*}
together with consistency, imply  the $8$ defining inequalities given by (\ref{eq:marg}) and (\ref{eq:joint}), as well as 
 $4$ non-negative inequalities of the form,
  \begin{equation}
    P(Y(x_0)=y,Y(x_1)=\tilde{y})\geq 0,
\end{equation}
where $y,\tilde{y} \in \{0,1\}$ and $4$ of the form,
 \begin{equation}\label{eq:obs-positivity}
 P(X(z_k)\!=\!i,Y(x_i)\!=\!y) = P(X\!=\!i,Y\!=\!y\;|\; Z\!=\!k) \geq 0,
 \end{equation}
 where $y,i\in \{0,1\}$;  the first equality in (\ref{eq:obs-positivity}) follows from (\ref{eq:z-obs-constraint}).
  
  \bigskip
  

\noindent\begin{tabular}{c@{ $\leq$ }c@{ \quad $\Leftrightarrow$\quad }rcl}
(a) & (A) & $P(Y\!=\!0, X\!=\!0\mid z_k)$& $\geq$ & $0$\\
(a) & (B) & $P(Y(x_1)\!=\!1)$ &$\leq$ & $1 - P(Y\!=\!0,X\!=\!1 \mid z_k)$\\
(a) & (C) & $P(Y(x_0)\!=\!0,Y(x_1)\!=\!0)$& $\geq$ & $0$  \\
(a) & (D) &$P(Y(x_0)\!=\!0,Y(x_1)\!=\!1)$& $\leq$ & $P(Y\!=\!0,X\!=\!0 \mid z_k) + P(Y\!=\!1,X\!=\!1 \mid z_k)$ \\[8pt]
(b) & (A) & $P(Y(x_1)\!=\!0)$& $\leq$ & $1 - P(Y\!=\!1,X\!=\!1 \mid z_k)$\\
(b) & (B) & $ P(Y\!=\!1,X\!=\!0 \mid z_k)$ &$\geq$ & $0$\\
(b) & (C) & $P(Y(x_0)\!=\!1,Y(x_1)\!=\!0)$& $\leq$ & $P(Y\!=\!1,X\!=\!0 \mid z_k) + P(Y\!=\!0,X\!=\!1 \mid z_k)$  \\
(b) & (D) & $P(Y(x_0)\!=\!1,Y(x_1)\!=\!1)$& $\geq$ & $0$ \\[8pt]
(c) & (A) &  $P(Y(x_0)\!=\!0,Y(x_1)\!=\!1)$ & $\geq$ & $0$\\
(c) & (B) & $P(Y(x_0)\!=\!1,Y(x_1)\!=\!1)$& $\leq$ & $P(Y\!=\!1,X\!=\!0 \mid z_k) + P(Y\!=\!1,X\!=\!1 \mid z_k)$\\
(c) & (C) & $P(Y(x_0)\!=\!1)$& $\leq$ & $1- P(Y\!=\!0,X\!=\!0 \mid z_k)$  \\
(c) & (D) & $P(Y\!=\!1,X\!=\!1 \mid z_k)$& $\geq$ & $0$ \\[8pt]
(d) & (A) & $P(Y(x_0)\!=\!0,Y(x_1)\!=\!0)$& $\leq$ & $P(Y\!=\!0,X\!=\!0 \mid z_k) + P(Y\!=\!0,X\!=\!0 \mid z_k)$\\
(d) & (B) & $P(Y(x_0)\!=\!1,Y(x_1)\!=\!0)$ & $\geq$ & $0$\\
(d) & (C) & $P(Y\!=\!0,X\!=\!1 \mid z_k)$& $\geq$ & $0$  \\
(d) & (D) &  $P(Y(x_0)\!=\!0)$& $\leq$ & $1- P(Y\!=\!1,X\!=\!0 \mid z_k)$ 
\end{tabular}


\section{Relation to the ACE$(X\rightarrow Y)$ bounds of Balke \& Pearl when $Z$ has two levels \label{app:two}}

We show explicitly here that in the case where $Z$ is binary, the bounds given by 
(\ref{eq:ace-bnds}) lead to those originally given by Balke and Pearl; see \citep{pearl:2000}.

We also show why it is the case that although the upper and lower bounds on
$P(Y(x_0)=1)$ and $P(Y(x_1) = 1)$ are each given by $4$ expressions, and --- as a 
consequence of variation independence --- these imply the bounds on the ACE, the latter
are each given by $8$ expressions, rather than $4\times 4 = 16$.

\subsection{Upper Bounds on $ACE (X\rightarrow Y)$}
With two levels of $Z$, the upper bounds on $P(Y(x_1)=1)$ are given by the four quantities defining $g(1,1)$:
\begin{align}
&1-P(x_1,y_0\mid z_0),\tag{$u_0^1$}\\ 
&1-P(x_1,y_0\mid z_1),\tag{$u_1^1$}\\ 
&P(y_1\mid z_0) + P(x_1,y_1\mid z_1) + P(x_0,y_0\mid z_1), \tag{$u_{01}^1$}\\ 
&P(y_1\mid z_1) + P(x_1,y_1\mid z_0) + P(x_0,y_0\mid z_0). \tag{$u_{10}^1$}  
\end{align}
The lower bounds on $P(Y(x_0)=1)$ are given by the four quantities corresponding to $1-g(0,0)$:
\begin{align}
&P(x_0,y_1\mid z_0),\tag{$l_0^0$}\\ 
&P(x_0,y_1\mid z_1),\tag{$l_1^0$}\\ 
&P(y_1\mid z_0) - P(x_1,y_1\mid z_1) - P(x_0,y_0\mid z_1), \tag{$l_{01}^0$}\\ 
&P(y_1\mid z_1) - P(x_1,y_1\mid z_0) - P(x_0,y_0\mid z_0). \tag{$l_{10}^0$} 
\end{align}

As noted earlier, owing to variation independence of $P(Y(x_0)=1)$ and $P(Y(x_1) = 1)$, we may compute upper bounds
on $ACE (X\rightarrow Y)$ simply by taking the upper bound on $P(Y(x_1) = 1)$ and subtracting the lower bound on $P(Y(x_0) = 1)$.

Table \ref{tab:two} relates the sixteen differences to the expressions given in (8.14b) on p.267 in  \citet{pearl:2000} and  in \citet{swanson:2018}, Table 3, row $A3+A4$.

\begin{table}[H]
\begin{tabular}{cccc}
Difference & Expression & Pearl & Swanson {\it et al.}\\[-6pt]
& or  why redundant &  {\scriptsize (8.14b) Row}  &  {\scriptsize Tab.3 A.3+A.4 Row}\\
\hline
$(u_0^1)-(l_0^0)$ & $P(x_0,y_0\mid z_0) + P(x_1,y_1\mid z_0)$ & \kern8pt 6\,%
\tablefootnote{In  \protect\citep{pearl:2000} this is given correctly as $p_{11.0} + p_{00.0}$, but in the
2009 reprinting there is a typo; it states $p_{11.0} - p_{00.0}$.}
& 1\\
$(u_0^1)-(l_1^0)$ & $1-P(x_1,y_0\mid z_0) - P(x_0,y_1\mid z_1)$ & 2 & \kern8pt 3\,%
\tablefootnote{In  \protect\citep{swanson:2018} Tab.3 A.3+A.4 Row 3 two different expressions are given and claimed to be equal.
The first, $p_{y_1|z_0} - p_{y_1|z_1} + p_{y_0,x_0|z_0} + p_{y_1,x_1|z_1}$ is correct, but the second contains a typo in a sign;
it should be the expression given here: $1-p_{y_0,x_1|z_0} \textcolor{blue}{-} p_{y_1,x_0|z_1}$.}
\\
$(u_0^1)-(l_{01}^0)$ & $ \left[(u_1^1)-(l_1^0)\right] + P(x_0,y_0\mid z_0)$ & -- & -- \\
$(u_0^1)-(l_{10}^0)$ & \multicolumn{1}{l}{$P(y_0\mid z_1) -P(x_1,y_0\mid z_0) $} & 3 & 5 \\[-4pt]
& \multicolumn{1}{r}{$ +\, P(x_1,y_1\mid z_0) + P(x_0,y_0\mid z_0)$} & & \\
\hline
$(u_1^1)-(l_0^0)$ & $1-P(x_1,y_0\mid z_1) - P(x_0,y_1\mid z_0)$ & 1 & \kern8pt 4\,%
\tablefootnote{In  \protect\citep{swanson:2018} Tab.3 A.3+A.4 Row 4 two different expressions are given and claimed to be equal.
The first, $p_{y_1|z_1} - p_{y_1|z_0} + p_{y_0,x_0|z_1} + p_{y_1,x_1|z_0}$ is correct, but the second contains a typo in a sign;
it should be the expression given here: $1-p_{y_0,x_1|z_1} \textcolor{blue}{-} p_{y_1,x_0|z_0}$.}
\\
$(u_1^1)-(l_1^0)$ & $P(x_0,y_0\mid z_1) + P(x_1,y_1\mid z_1)$ & 5 & 2\\
$(u_1^1)-(l_{01}^0)$ & \multicolumn{1}{l}{$P(y_0\mid z_0) -P(x_1,y_0\mid z_1) $} & 4 & 6 \\[-4pt]
& \multicolumn{1}{r}{$ +\, P(x_1,y_1\mid z_1) + P(x_0,y_0\mid z_1)$} & & \\
$(u_1^1)-(l_{10}^0)$ &   $\left[(u_0^1)-(l_0^0)\right] + P(x_0,y_0\mid z_0)$  & -- & -- \\
\hline 
$(u_{01}^1)-(l_0^0)$ &  $\left[(u_1^1)-(l_1^0)\right] + P(x_1,y_1\mid z_0)$ & -- & -- \\
$(u_{01}^1)-(l_1^0)$ &\multicolumn{1}{l}{$P(y_1\mid z_0) - P(x_0,y_1\mid z_1)$}& 7 & 7 \\[-4pt]
& \multicolumn{1}{r}{$ +\, P(x_1,y_1\mid z_1) + P(x_0,y_0\mid z_1)$}\\
$(u_{01}^1)-(l_{01}^0)$ & $2\!\left[ (u_1^1)-(l_1^0) \right]$ & -- & -- \\
$(u_{01}^1)-(l_{10}^0)$ & \multicolumn{1}{l}{$\left[(u_0^1)-(l_1^0)\right]$}  &-- & --   \\[-4pt]
&  \multicolumn{1}{r}{$+ P(x_1,y_1\mid z_0) + P(x_0,y_0\mid z_1)$} & &\\
\hline
$(u_{10}^1)-(l_0^0)$ &\multicolumn{1}{l}{$P(y_1\mid z_1) - P(x_0,y_1\mid z_0)$}& 8 & 8 \\[-4pt]
& \multicolumn{1}{r}{$ +\, P(x_1,y_1\mid z_0) + P(x_0,y_0\mid z_0)$}\\
$(u_{10}^1)-(l_1^0)$ &  $\left[(u_0^1)-(l_0^0)\right] + P(x_1,y_1\mid z_1)$ & -- & -- \\
$(u_{10}^1)-(l_{01}^0)$ & \multicolumn{1}{l}{$\left[(u_1^1)-(l_0^0)\right]$}  &-- & --   \\[-4pt]
&  \multicolumn{1}{r}{$+ P(x_1,y_1\mid z_1) + P(x_0,y_0\mid z_0)$} & &\\
$(u_{10}^1)-(l_{10}^0)$ & $2\!\left[ (u_0^1)-(l_0^0) \right]$ & -- & -- \\
\end{tabular}
\caption{The upper bounds on ACE$(X\rightarrow Y)$ in the case where $Z$ is binary, that arise by taking the differences between the four upper bounds on  $P(Y(x_1) = 1)$ and the four lower bounds on $P(Y(x_0) = 1)$; an explanation is given as to why $8$ differences are redundant. See footnotes on the next page. \label{tab:two}}
\end{table}


\subsection{Lower Bounds on $ACE (X\rightarrow Y)$}
With two levels of $Z$, the lower bounds on $P(Y(x_1)=1)$ are given by the four quantities defining $1-g(1,0)$:
\begin{align}
&P(x_1,y_1\mid z_0),\tag{$l_0^1$}\\ 
&P(x_1,y_1\mid z_1),\tag{$l_1^1$}\\ 
&P(y_1\mid z_0) - P(x_1,y_0\mid z_1) - P(x_0,y_1\mid z_1), \tag{$l_{01}^1$}\\ 
&P(y_1\mid z_1) - P(x_1,y_0\mid z_0) - P(x_0,y_1\mid z_0). \tag{$l_{10}^1$}  
\end{align}
The upper bounds on $P(Y(x_0)=1)$ are given by the four quantities corresponding to $g(0,1)$:
\begin{align}
&1-P(x_0,y_0\mid z_0),\tag{$u_0^0$}\\ 
&1-P(x_0,y_0\mid z_1),\tag{$u_1^0$}\\ 
&P(x_0,y_1\mid z_0) + P(x_1,y_0\mid z_0) + P(y_1\mid z_1), \tag{$u_{01}^0$}\\ 
&P(x_0,y_1\mid z_1) + P(x_1,y_0\mid z_1) + P(y_1\mid z_0). \tag{$u_{10}^0$} 
\end{align}

Again, owing to variation independence of $P(Y(x_0)=1)$ and $P(Y(x_1) = 1)$, we may compute lower bounds
on $ACE (X\rightarrow Y)$ simply by taking the lower bound on $P(Y(x_1) = 1)$ and subtracting the upper bound on $P(Y(x_0) = 1)$.

\begin{table}[H]
\begin{tabular}{cccc}
Difference & Expression & Pearl & Swanson {\it et al.}\\[-6pt]
& or  why redundant &  {\scriptsize (8.14a) Row}  &  {\scriptsize Tab.2 A.3+A.4 Row}\\
\hline
$(l_0^1)-(u_0^0)$ & $-P(x_1,y_0\mid z_0) - P(x_0,y_1\mid z_0)$ & 6 & 1\\
$(l_0^1)-(u_1^0)$ &  $P(x_1,y_1\mid z_0) + P(x_0,y_0\mid z_1) -1 $ & 2 & 3\\
$(l_0^1)-(u_{01}^0)$ & \multicolumn{1}{l}{$-P(y_1\mid z_1) + P(x_1,y_1\mid z_0) $} & 3 & 5 \\[-4pt]
& \multicolumn{1}{r}{$ -\, P(x_0,y_1\mid z_0) - P(x_1,y_0\mid z_0)$} & & \\
$(l_0^1)-(u_{10}^0)$ & $ \left[(l_1^1)-(u_1^0)\right] - P(x_0,y_1\mid z_0)$ & -- & -- \\
\hline
$(l_1^1)-(u_0^0)$ & $P(x_1,y_1\mid z_1) + P(x_0,y_0\mid z_0) -1$ & 1 & 4\\
$(l_1^1)-(u_1^0)$ & $-P(x_1,y_0\mid z_1) - P(x_0,y_1\mid z_1)$  & 5 & 2\\
$(l_1^1)-(u_{01}^0)$ & $ \left[(l_0^1)-(u_0^0)\right] - P(x_0,y_1\mid z_1)$ & -- & -- \\
$(l_1^1)-(u_{01}^0)$ & \multicolumn{1}{l}{$-P(y_1\mid z_0) + P(x_1,y_1\mid z_1) $} & 4 & 6 \\[-4pt]
& \multicolumn{1}{r}{$ -\, P(x_0,y_1\mid z_1) - P(x_1,y_0\mid z_1)$} & & \\
\hline
$(l_{01}^1)-(u_0^0)$ &  $\left[(l_1^1)-(u_1^0)\right] - P(x_1,y_0\mid z_0)$ & -- & -- \\
$(l_{01}^1)-(u_1^0)$ &\multicolumn{1}{l}{$-P(y_0\mid z_0) + P(x_0,y_0\mid z_1)$}& 7 & 7 \\[-4pt]
& \multicolumn{1}{r}{$ -\, P(x_1,y_0\mid z_1) - P(x_0,y_1\mid z_1)$}\\
$(l_{01}^1)-(u_{01}^0)$ & \multicolumn{1}{l}{$\left[(l_0^1)-(u_1^0)\right]$}  &-- & --   \\[-4pt]
&  \multicolumn{1}{r}{$- P(x_1,y_0\mid z_0) - P(x_0,y_1\mid z_1)$} & &\\
$(l_{01}^1)-(u_{10}^0)$ & $2\!\left[ (l_1^1)-(u_1^0) \right]$ & -- & -- \\
\hline
$(l_{10}^1)-(u_0^0)$ &\multicolumn{1}{l}{$-P(y_0\mid z_1) + P(x_0,y_0\mid z_0)$}& 8 & 8 \\[-4pt]
& \multicolumn{1}{r}{$ -\, P(x_1,y_0\mid z_0) - P(x_0,y_1\mid z_0)$}\\
$(l_{10}^1)-(u_1^0)$ &   $\left[(l_0^1)-(u_0^0)\right] - P(x_1,y_0\mid z_1)$ & -- & -- \\
$(l_{10}^1)-(u_{01}^0)$ & $2\!\left[ (l_0^1)-(u_0^0) \right]$ & -- & -- \\
$(l_{10}^1)-(u_{10}^0)$ &  \multicolumn{1}{l}{$\left[(l_1^1)-(u_0^0)\right]$}  &-- & --   \\[-4pt]
&  \multicolumn{1}{r}{$ - P(x_0,y_1\mid z_0) - P(x_1,y_0\mid z_1)$} & &\\
\hline
\end{tabular}
\caption{The lower bounds on ACE$(X\rightarrow Y)$ arising by taking the differences between the four lower bounds on  $P(Y(x_1) = 1)$ and the four upper bounds on $P(Y(x_0) = 1)$; an explanation is given as to why $8$ differences are redundant. \label{tab:three}}
\end{table}

 Table \ref{tab:three}  relates the sixteen differences to the expressions given in (8.14a) on p.267 in  \citet{pearl:2009, pearl:2000} and \citet{swanson:2018}, Table 2, row $A3+A4$.

\end{document}